\documentclass[sn-mathphys]{sn-jnl}% Math and Physical Sciences Reference Style
\usepackage{graphicx}%
\usepackage{multirow}%
\usepackage{amsmath,amssymb,amsfonts}%
\usepackage{amsthm}%
\usepackage{mathrsfs}%
\usepackage[title]{appendix}%
\usepackage{xcolor}%
\usepackage{textcomp}%
\usepackage{manyfoot}%
\usepackage{booktabs}%
\usepackage{algorithm}%
\usepackage{listings}%
\begin{document}

\title[automatic projection parameter increase]{Achieving binary topology optimization solutions via automatic projection parameter increase}


\author*[1]{\fnm{Peter} \sur{Dunning}}\email{peter.dunning@abdn.ac.uk}

\affil[1]{\orgdiv{School of Engineering}, \orgname{University of Aberdeen}, \orgaddress{\street{King's College}, \city{Aberdeen}, \postcode{AB24 3UE}, \country{UK}}}

\abstract{A method is created to automatically increase the threshold projection parameter in three-field density-based topology optimization to achieve a near binary design. The parameter increase each iteration is based on an exponential growth function, where the growth rate is dynamically changed during optimization by linking it to the change in objective function. This results in a method that does not need to be tuned for specific problems, or optimizers, and the same set of hyper-parameters can be used for a wide range of problems. The effectiveness of the method is demonstrated on several 2D benchmark problems, including linear buckling and geometrically nonlinear problems.}

\keywords{Topology optimization, Threshold projection, Linear buckling, Nonlinear geometry}

\maketitle

\section{Introduction}\label{sec1}

Topology optimization aims to find the ideal layout, shape and size of a structure for given objectives and constraints. These problems are typically framed by deciding where material should exist within the design domain, leading to a binary value problem, where a value of 1 indicates material, and 0 is empty space (or void). However, solving binary optimization problems is challenging, and it is therefore popular to relax the problem using continuous design variables, as introduced by \cite{bendsoe1988generating}. The most popular topology optimization approach using this principle is the Simple Isotropic Material with Penalization (SIMP) method, where a pseudo-density value is assigned to each element, which can vary between 0 (void) and 1 (solid) \citep{bendsoe1989optimal, rozvany1992generalized}.

The pseudo-density method is usually coupled with filtering techniques to avoid numerical issues such as checkerboard patterns and mesh dependency \citep{sigmund1998numerical}. A popular approach is to use a density filter, as this results in a consistent problem formulation, where the derivatives of the optimization objectives and constraints can be explicitly computed \citep{bruns2001topology,bourdin2001filters}. However, using a density filter can lead to solutions with a significant amount intermediate pseudo-density values (between 0 and 1), far from the ideal binary design \citep{wang2011projection}. Intermediate, or "gray" elements can also appear in problems where the optimizer may exploit the properties of these elements to improve the objective, or meet a constraint. For example, intermediate material may appear to create a low-density core to improve buckling performance, similar to a sandwich structure.

To force a design to be binary, a further operation can be used to threshold pseudo-density values such that any value above the threshold becomes 1 (solid) and any below becomes 0 (void) \citep{guest2004achieving}. However, to use a threshold during gradient-based optimization, the threshold projection function is also relaxed to take a continuous form \citep{guest2004achieving,xu2010volume,wang2011projection}. These relaxed formulations use a projection parameter that defines the sharpness of the projection, where a small value results in little, or no change to pseudo-density values, and a value of infinity results in a binary design (albeit non-differentiable). This introduces an optimization hyper-parameter that must be chosen by the user. From experience, many researchers use a continuation scheme for the projection parameter, where it is kept at a small value (e.g. 1) during early  optimization iterations to allow a design to emerge. It is then incrementally increased to push the design towards a binary solution, for example see \citep{ferrari2020new, dunning2023stability}. These continuation schemes are usually ad-hoc and heuristic, needing to be tuned for each new problem. It is also not clear if these heuristic schemes achieve the best balance between computational cost and performance of the design (including how discrete it is).

An alternative approach was introduced by \cite{da2023topology} where the threshold projection parameter is included as an extra element-wise design variable. The lower limit for each variable is dynamically updated based on the current grayness measure (a measure on non-discreetness) for the element. In addition, material penalisation parameters are treated in the same way. The method is an improvement on heuristic continuation schemes, and works well for the problem studied (fibre reinforced structures with stress constraints), but still has several tuning parameters and it is unclear if the method would work on a wider class of problems.

Another idea was recently proposed by \cite{ha2024automatic}  where a surrogate optimization method is used as an alternative to manually determining topology optimization hyper-parameters, including continuation parameters. This involves solving the original topology optimization problem as an inner problem in a nested loop, where hyper-parameters are optimized in an outer loop aided by machine learning techniques. This can reduce the effort required to manually tune hyper-parameters, but may not be feasible for computationally expensive problems (e.g. nonlinear and multi-disciplinary problems).

Note that several topology optimization methods have been developed that use binary design variables, such as Bi-directional Evolutionary Structural Optimization (BESO) \citep{huang2007convergent}, and more recently Topology Optimization of Binary Structures (TOBS) \citep{sivapuram2018topology}, amongst others. However, these methods have their own heuristic hyper-parameters that may be problem dependent. There are also some methods that allow intermediate densities in the solution by linking their properties to a porous micro-structure, as originally explored in the seminal paper by \cite{bendsoe1988generating}. This idea has been recently applied to develop multi-scale methods, see for example \citep{christensen2023topology,hubner2023two}.

In this paper, a simple alternative method for obtaining near binary topology optimization designs is proposed. Starting with the three-field density-based method (where a density filter is combined with threshold projection) a method is developed to automatically increase the threshold projection parameter during optimization by linking it to the progress of the objective function. The three-field density-based method is reviewed in Section \ref{sec21} and the new automatic continuation scheme is introduced in Sections \ref{sec22} and \ref{sec23}. The method is demonstrated using several benchmark problems in Section \ref{sec3}, including comparison with heuristic, ad-hoc continuation schemes from the literature. 

\section{Methodology}\label{sec2}

\subsection{Three-field topology optimization}\label{sec21}
The three-field density-based topology optimization method starts with the design variable field, where a scalar value $x_e \in [ 0, 1 ]$ is assigned to each element, $e$, in the mesh. A density-based filter is then applied to avoid well-known issues, such as checkerboard patterns and mesh dependent solutions \citep{bruns2001topology,bourdin2001filters}. 
\begin{equation}
\label{filter}
\tilde{x}_e =  \frac{\sum_{i \in N_{i,e}}{H_{i,e} x_i}}{\sum_{i \in N_{i,e}}{H_{i,e}}}
\end{equation}
where $N_{i,e}$ is the set of elements $i$ where the distance between the center of element $i$ and center of element $e$, $d(e,i)$, is less than the filter radius $r_{\min}$. $H_{i,e}$  is a weight factor, defined as:
\begin{equation}
\label{H_filt}
H_{i,e} = \max (0, r_{\min} - d(e,i) )
\end{equation}
The filtered field $\tilde{\boldsymbol{x}}$ creates a band of intermediate density values around the boundary of the structure. Thus, a threshold projection function is used to obtain a more crisp (or binary) design description \citep{xu2010volume,wang2011projection}:
\begin{equation}
\label{H_proj}
\bar{x}_e = \frac{\tanh{(\beta \eta)} + \tanh{(\beta (\tilde{x}_e - \eta))}}{\tanh{(\beta \eta)} + \tanh{(\beta (1-\eta))}}
\end{equation}
where $\bar{x}$ is the physical density field, $\beta$ is the projection parameter that controls the sharpness of the projection, and $\eta$ is the threshold value, typically set as $\eta = 0.5$.

The value of $\bar{x}_e$ is used to compute the physical properties of the element, where the volume is directly proportional to $\bar{x}_e$ and the modified SIMP method is used to compute the Young's Modulus:
\begin{equation}
\label{SIMP}
E_e (\bar{x}_i) =  E_{\min} + \bar{x}_e^p (E_0 - E_{\min})
\end{equation}
where $E_0$ is the Young's modulus of the material, $E_{\min}$ is a small value to prevent the stiffness matrix becoming singular and $p$ is the penalization factor, where the typical value of $p=3$ is used in this paper. Note that the physical pseudo-density values $\bar{x}_e$ are used to plot solutions.

\subsection{Binary solutions via threshold projection}\label{sec22}
The solution to a topology optimization problem should ideally be binary, indicating where material exists ($\bar{x}_e = 1$), and where is does not ($\bar{x}_e = 0$). Therefore, we aim to avoid intermediate density values in the solution, although they are necessary during optimization to make the problem differential such that it can be solved by efficient gradient-based optimizers. 

The gray level indicator function \citep{sigmund2007morphology} can be used to measure how close a design is to being binary:
\begin{equation}
\label{gray}
G(\bar{x}) =  \frac{4}{N} \sum_{e=1}^{N}{\bar{x}_e (1 - \bar{x}_e)}
\end{equation}
where $N$ is the number of elements. Smaller values of $G(\bar{x})$ indicate the design is closer to being completely binary, with $G(\bar{x}) = 0$ characterizing a true binary design.

The threshold projection function, Eq. (\ref{H_proj}), can be used to reduce the gray level of a design by increasing the projection parameter $\beta$, thus making it closer to being binary, regardless of other parameters, such as the penalization power $p$ in Eq. (\ref{SIMP}). When $\beta$ is increased, then if $\eta < \bar{x}_e < 1$ it increases to be closer to 1, and if $0 < \bar{x}_e < \eta$ it decreases to become closer to 0. In the limit of $\beta = \infty$, only three physical density values are possible: $\bar{x}_e \in \{ 0, 0.5, 1 \}$, where $\bar{x}_e = 0.5$ is only possible when $\tilde{x}_e = \eta$, thus the design is practically binary.

However, starting the optimization with a large value of $\beta$ can cause convergence problems and limits design exploration. Therefore, researchers often start with a small value of $\beta$ and employ continuation schemes to gradually increase it during optimization. This allows a design to emerge in the early iterations, whilst being projected close to a binary design at the end. Current continuation approaches for the projection parameter are dependent on several hyper-parameters, which are often tuned for each new problem because they can significantly influence the performance of the solution and number of iterations. This issue is addressed in the following section.

\subsection{Automatic projection parameter increase}\label{sec23}
Some continuation schemes increase $\beta$ after a number of iteration steps. However, this can cause a discontinuity in the optimization process. There are also questions on how frequent $\beta$ should be increased, by how much should it be increased, and when to stop increasing. The method proposed here addresses the question of how frequent, by checking for $\beta$ increase every iteration. The amount of increase, $\Delta \beta$, is linked to the rate of decrease in the objective function, as this indicates whether a design is still emerging (large relative change in objective), or is settling (small relative change in objective). Assuming the objective is to minimize a function, $f$, the increase in $\beta$ at iteration $k$ is calculated by:
\begin{equation}
\label{delbeta}
\Delta \beta^k =  \max \left(-\frac{\gamma}{2}~\frac{f^k + f^{k-1}}{f^k - f^{k-1}} , 0 \right )
\end{equation}
where $\gamma$ is a tuning parameter that sets how fast $\beta$ increases. A smaller value leads to a slower increase, allowing the optimizer more time to refine the design, before being pushed towards a binary solution, but potentially taking more iterations. Whereas a bigger value pushes the design towards a binary solution quicker, at the possible sacrifice in objective function value. The value of $\beta$ in the next iteration is then: 
\begin{equation}
\label{newbeta}
\beta^{k+1} =  \beta^k + \min \left ( \Delta \beta^k , \Delta \beta_{max} \right )
\end{equation}
where $\Delta \beta_{max}$ is the maximum possible increase. Note that a maximum value on $\beta$ should also be set to avoid numerical issues with extremely large numbers, and that $\beta=1.0$ in the first iteration.

Although the proposed automatic increase method is not free from tuning parameters, testing has shown that setting $\gamma = 10^{-4}$ and $\Delta \beta_{max}=0.2\beta^k$ works well for a variety of problems, as shown in Section \ref{sec3} below. Therefore, the need to tune these parameters for new problems is reduced, compared to existing methods, and in most cases tuning is not required.

The final part of the proposed method is to use the current value of the gray level indicator function, Eq. (\ref{gray}), to determine if $\beta$ should be increased and if the optimization process can be stopped. Thus, $\beta$ is only increased if, $G(\bar{x}) > \epsilon$, and the optimization process should continue until, $G(\bar{x}) \leq \epsilon$, and all other stopping, or convergence criteria are met. The value of $\epsilon$ can be set by the user to trade-off computational time (bigger $\epsilon$ leads to fewer iterations) with how close the solution is to a binary design (smaller $\epsilon$ means less gray in the solution). For all examples in this work $\epsilon = 0.01$, providing a reasonable balance between computational cost and grayness.

\section{Examples}\label{sec3}

\begin{figure}[h]
    \centering
    \includegraphics[width=0.8\textwidth]{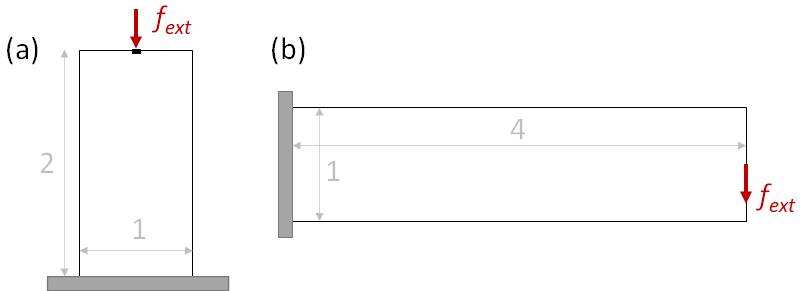}
    \caption{Example problems, a) compressed column, b) cantilever.}
    \label{examples}
\end{figure}

Several benchmark problems are used to demonstrate the method described in Section \ref{sec23}. The first is a compressed column, Figure \ref{examples}a, and the second is a cantilever with aspect ratio 4, Figure \ref{examples}b.

The compressed column is discretized using $120 \times 240$ four-node plane stress elements with a thickness of 1. The applied load is $10^{-3}$, which is spread over 8 elements, $E_0 = 1$, $E_{min} = 10^{-6}$ and Poisson's ratio $\nu = 0.3$. In addition, a small region around the applied load ($8 \times 4$ elements) is fixed to remain part of the structure by setting $\bar{x}_e = 1$. Two different filter radii in Eq. (\ref{H_filt}) (4, or 8 times element edge length, $h$) are used to test the ability of the proposed method to achieve near binary designs with different length scales, noting that a larger filter radius creates more gray elements around the boundary.

The cantilever is discretized using three mesh sizes: $160 \times 40$, $320 \times 80$, and $640 \times 160$ to test the ability of the proposed method to achieve a near binary design whilst refining the mesh. This problem is solved considering geometrically nonlinear behaviour using the corotational element proposed by \cite{crisfield1996co}. Two different applied loads are used: $2 \times 10^{5}$ and $3 \times 10^{5}$. The element thickness is 0.1, $E_0 = 3 \times 10^9$, $E_{min} = 3$ and $\nu = 0.4$. The filter radius in Eq. (\ref{H_filt}) is set to 0.075.

Note that the same stopping criteria are used for all problems, where the gray level indicator function, Eq. (\ref{gray}), must be below 0.01, the absolute value of relative change in objective function between two iterations must be below $10^{-5}$, and all constraints are satisfied.

\subsection{Compressed column results}\label{sec31}
The compressed column example is solved for two different problems. The first is to maximize the buckling load factor (or minimize the inverse) with a 35\% volume constraint, using a uniform density starting design equal to the volume constraint. The second is to minimize volume with a buckling load factor constraint (minimum value of 15) and compliance constraint - set at 2 times the compliance for the initial fully solid design. Both problems are solved using the Matlab code developed by \cite{ferrari2021topology}, which includes an inbuilt continuation method for incrementally increasing $\beta$ during optimization. For comparison, the problems are also solved using this in-built continuation method. The default settings are that $\beta = 1$ for the first 400 iterations, it is then increased by 2 every 25 iterations until a maximum value of 25. A further comparison is made to a ``hand tuned" modified version of this approach following analysis of the optimization history using the default settings, where initially $\beta = 1$ for the first 200 iterations, it is then increased by 2 every 25 iterations up to a possible maximum value of 500. The problems are solved using the inbuilt optimally criteria method in the Matlab code \citep{ferrari2021topology}. To avoid issues with mode switching, 30 modes are considered for the maximization of buckling load factor problem, and 20 for the buckling constraint in the volume minimization problem. In both cases, the maximum buckling load factor is approximated using a K-S function with an aggregation parameter of 160. 

\begin{figure}[h]
    \centering
    \includegraphics[width=0.75\textwidth]{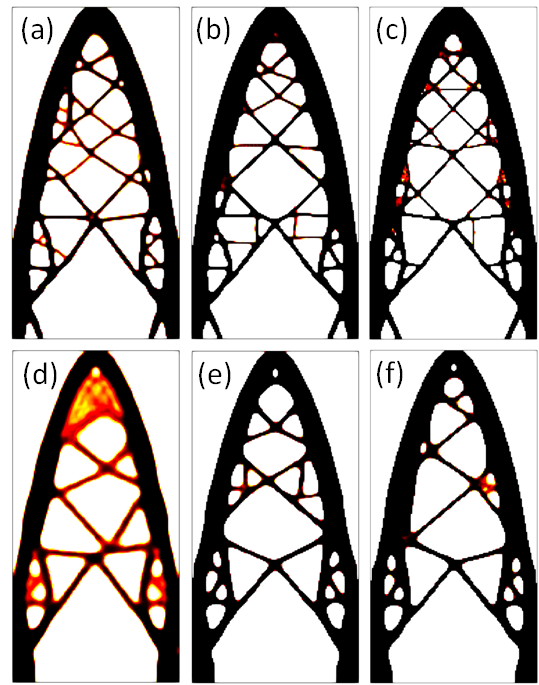}
    \caption{Compressed column results for maximization of buckling load factor. Solutions with filter radius = $4h$: a) Default continuation scheme, b) modified scheme, c) automatic scheme. Solutions with filter radius = $8h$: d) Default continuation scheme, e) modified scheme, f) automatic scheme.}
    \label{BVsol}
\end{figure}

The solutions for the maximum buckling load factor problem are shown in Figure \ref{BVsol}. Note that the solutions for the default continuation approach are stopped after 2000 iterations, as the maximum value of $\beta = 25$ could not meet the gray level stopping criteria - i.e. the default maximum value of $\beta$ is too small to achieve a near binary design. This is most noticeable in the solution using a filter radius of $8h$ in Figure \ref{BVsol}d, where large areas of intermediate density material exist, demonstrating the challenge in achieving near binary designs for this problem. The modified scheme allows for a bigger $\beta$ values and a solution for the filter radius $8h$ problem that meets the gray level criterion is obtained after 2751 iterations, Figure \ref{BVsol}e. The proposed automatic scheme also produces a design that meets the gray level criterion, Figure \ref{BVsol}f, but it only needs 1029 iterations. The designs are slightly different, but the maximum buckling load factor for the design obtained using the automatic scheme is slightly better (16.26 versus 16.22 - see Table \ref{BVtab}). Similar results can be seen when the filter radius is $4h$, although the reduction in number of iterations comparing the automatic and modified schemes is smaller, as the smaller filter radius makes it easier to meet the gray level stopping criterion. This suggests that the modified scheme could be further improved for the $8h$ filter radius problem, but would require further optimization runs and hand tuning, which are not required for the proposed automatic scheme.

\begin{figure}[h]
    \centering
    \includegraphics[width=0.95\textwidth]{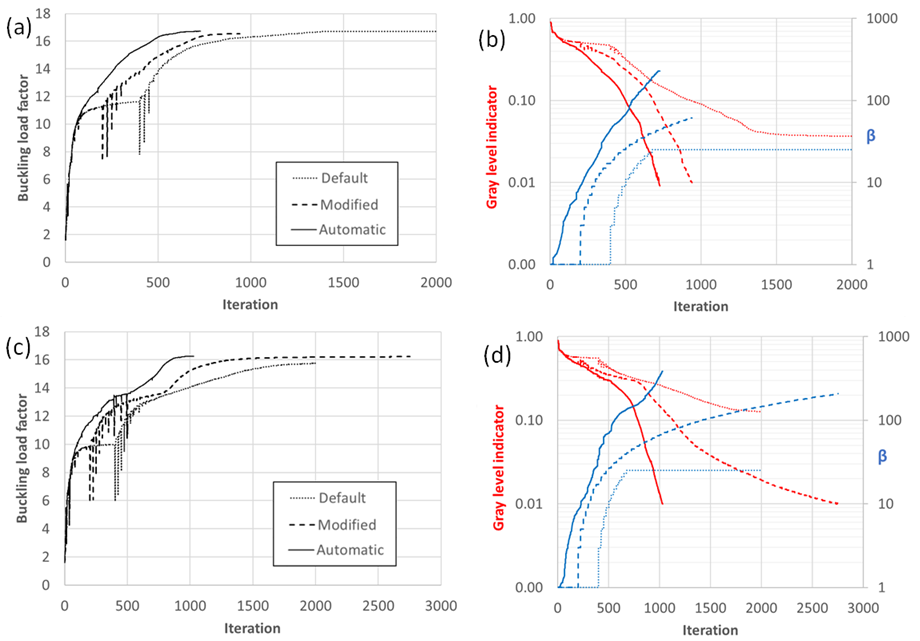}
    \caption{Compressed column convergence history for maximization of buckling load factor. Convergence of volume objective for: a) filter radius = $4h$, c) filter radius = $8h$. Convergence of gray level indicator function, $G(x)$ , and threshold projection parameter, $\beta$, for b) filter radius = $4h$, d) filter radius = $8h$.}
    \label{BVplot}
\end{figure}

Optimization histories for all problems, including the progress of the the gray level function and $\beta$ value, are shown in Figure \ref{BVplot}. They demonstrate how the automatic scheme modifies $\beta$ in relation to the progress of the objective function, and is able to push the design towards a near binary solution (as shown by the gray level indicator) quickly, without compromising the performance of the solution, in comparison with the other schemes.

\begin{table}[h]
\begin{tabular}{lcccc}
\textbf{Filter radius = $4h$}             & \multicolumn{1}{l}{}                      & \multicolumn{1}{l}{}            & \multicolumn{1}{l}{}         & \multicolumn{1}{l}{}                      \\ \hline
\multicolumn{1}{|l|}{Continuation method} & \multicolumn{1}{c|}{Buckling load factor} & \multicolumn{1}{c|}{Iterations} & \multicolumn{1}{c|}{$\beta$} & \multicolumn{1}{c|}{Gray value}           \\ \hline
\multicolumn{1}{|l|}{Default}            & \multicolumn{1}{c|}{16.72}                & \multicolumn{1}{c|}{2000}       & \multicolumn{1}{c|}{25}      & \multicolumn{1}{c|}{$3.68 \time 10^{-2}$} \\ \hline
\multicolumn{1}{|l|}{Modified}            & \multicolumn{1}{c|}{16.54}                & \multicolumn{1}{c|}{941}        & \multicolumn{1}{c|}{61}      & \multicolumn{1}{c|}{$1.00 \time 10^{-2}$} \\ \hline
\multicolumn{1}{|l|}{Automatic}           & \multicolumn{1}{c|}{16.72}                & \multicolumn{1}{c|}{728}        & \multicolumn{1}{c|}{229.7}   & \multicolumn{1}{c|}{$9.10 \time 10^{-3}$} \\ \hline
                                          & \multicolumn{1}{l}{}                      & \multicolumn{1}{l}{}            & \multicolumn{1}{l}{}         & \multicolumn{1}{l}{}                      \\
\textbf{Filter radius = $8h$}             & \multicolumn{1}{l}{}                      & \multicolumn{1}{l}{}            & \multicolumn{1}{l}{}         & \multicolumn{1}{l}{}                      \\ \hline
\multicolumn{1}{|l|}{Continuation method} & \multicolumn{1}{c|}{Buckling load factor} & \multicolumn{1}{c|}{Iterations} & \multicolumn{1}{c|}{$\beta$} & \multicolumn{1}{c|}{Gray value}           \\ \hline
\multicolumn{1}{|l|}{Default}     & \multicolumn{1}{c|}{15.75}                & \multicolumn{1}{c|}{2000}       & \multicolumn{1}{c|}{25}      & \multicolumn{1}{c|}{$1.27 \time 10^{-1}$} \\ \hline
\multicolumn{1}{|l|}{Modified}            & \multicolumn{1}{c|}{16.22}                & \multicolumn{1}{c|}{2751}       & \multicolumn{1}{c|}{207}     & \multicolumn{1}{c|}{$9.87 \time 10^{-3}$} \\ \hline
\multicolumn{1}{|l|}{Automatic}           & \multicolumn{1}{c|}{16.26}                & \multicolumn{1}{c|}{1029}       & \multicolumn{1}{c|}{385.9}   & \multicolumn{1}{c|}{$9.93 \time 10^{-3}$} \\ \hline
\end{tabular}
\caption{Results for compressed column for maximization of buckling load factor.}
\label{BVtab}
\end{table}

\begin{figure}[h]
    \centering
    \includegraphics[width=0.75\textwidth]{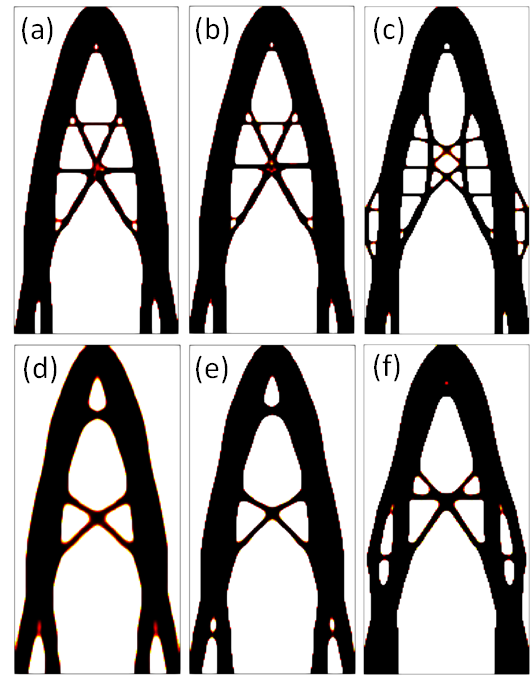}
    \caption{Compressed column results for minimization of volume. Solutions with filter radius = $4h$: a) Default continuation scheme, b) modified scheme, c) automatic scheme. Solutions with filter radius = $8h$: d) Default continuation scheme, e) modified scheme, f) automatic scheme.}
    \label{VCBsol}
\end{figure}

\begin{figure}[h]
    \centering
    \includegraphics[width=0.95\textwidth]{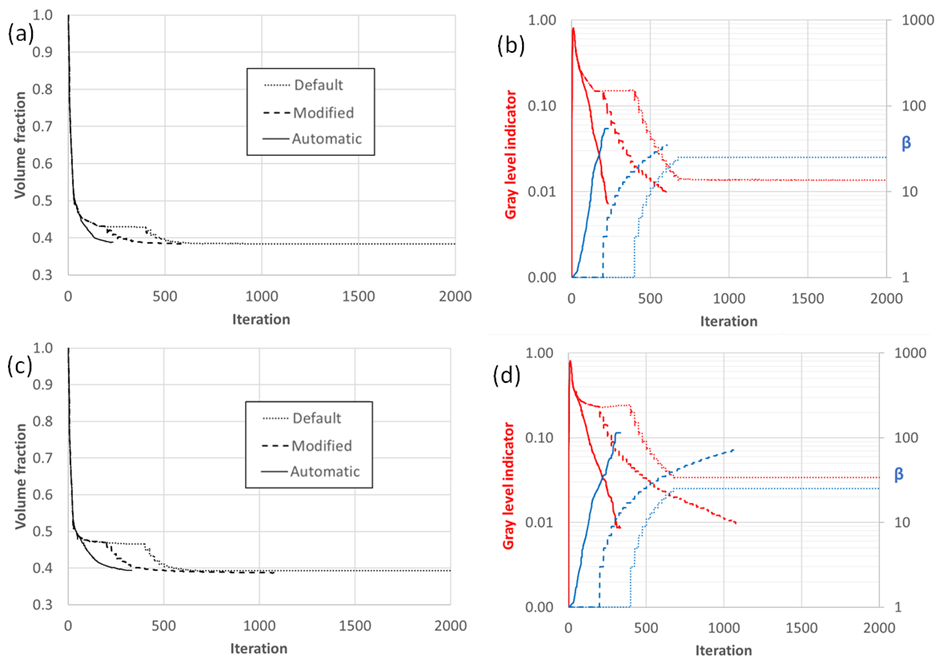}
    \caption{Compressed column convergence history for minimization of volume. Convergence of volume objective for: a) filter radius = $4h$, c) filter radius = $8h$. Convergence of gray level indicator function, $G(x)$, and threshold projection parameter, $\beta$, for b) filter radius = $4h$, d) filter radius = $8h$.}
    \label{VCBplot}
\end{figure}

The solutions of the volume minimization problem are shown in Figure \ref{VCBsol}. Again, the default scheme is unable to produce designs that meet the gray level stopping criteria, and are stopped after 2000 iterations. The modified and automatic schemes do produce designs that meet the criterion, although the automatic scheme achieves this in significantly fewer iterations (approximately a third - see Table \ref{VCBtab}). However, for this examples, the slower, modified scheme obtains slightly better solutions with smaller volume fractions (by approximately 1-2\%). This is evident from the simpler designs obtained by the modified scheme (e.g. compare Figure \ref{VCBsol}b with c, or Figure \ref{VCBsol}e with f). However, this small increase in performance, comes at significantly more computational cost, as seen by the required number of iterations in Table \ref{VCBtab}. The optimization histories, Figure \ref{VCBplot}, show the same as the previous problem, where the automatic scheme is able to adapt to the optimization progress and push the design towards a near binary solution quicker than the other two schemes, which rely on a predetermined $\beta$ increase scheme.

\begin{table}[h]
\begin{tabular}{lcccc}
\textbf{Filter radius = $4h$}             & \multicolumn{1}{l}{}                 & \multicolumn{1}{l}{}            & \multicolumn{1}{l}{}         & \multicolumn{1}{l}{}                      \\ \hline
\multicolumn{1}{|l|}{Continuation method} & \multicolumn{1}{c|}{Volume fraction} & \multicolumn{1}{c|}{Iterations} & \multicolumn{1}{c|}{$\beta$} & \multicolumn{1}{c|}{Gray value}           \\ \hline
\multicolumn{1}{|l|}{Default}            & \multicolumn{1}{c|}{0.384}           & \multicolumn{1}{c|}{2000}       & \multicolumn{1}{c|}{25}      & \multicolumn{1}{c|}{$1.36 \time 10^{-2}$} \\ \hline
\multicolumn{1}{|l|}{Modified}            & \multicolumn{1}{c|}{0.384}           & \multicolumn{1}{c|}{606}        & \multicolumn{1}{c|}{35}      & \multicolumn{1}{c|}{$9.45 \time 10^{-3}$} \\ \hline
\multicolumn{1}{|l|}{Automatic}           & \multicolumn{1}{c|}{0.388}           & \multicolumn{1}{c|}{232}        & \multicolumn{1}{c|}{54.5}    & \multicolumn{1}{c|}{$7.29 \time 10^{-3}$} \\ \hline
                                          & \multicolumn{1}{l}{}                 & \multicolumn{1}{l}{}            & \multicolumn{1}{l}{}         & \multicolumn{1}{l}{}                      \\
\textbf{Filter radius = $8h$}             & \multicolumn{1}{l}{}                 & \multicolumn{1}{l}{}            & \multicolumn{1}{l}{}         & \multicolumn{1}{l}{}                      \\ \hline
\multicolumn{1}{|l|}{Continuation method} & \multicolumn{1}{c|}{Volume fraction} & \multicolumn{1}{c|}{Iterations} & \multicolumn{1}{c|}{$\beta$} & \multicolumn{1}{c|}{Gray value}           \\ \hline
\multicolumn{1}{|l|}{Default}            & \multicolumn{1}{c|}{0.393}           & \multicolumn{1}{c|}{2000}       & \multicolumn{1}{c|}{25}      & \multicolumn{1}{c|}{$3.39 \time 10^{-2}$} \\ \hline
\multicolumn{1}{|l|}{Modified}            & \multicolumn{1}{c|}{0.387}           & \multicolumn{1}{c|}{1076}       & \multicolumn{1}{c|}{73}      & \multicolumn{1}{c|}{$9.71 \time 10^{-3}$} \\ \hline
\multicolumn{1}{|l|}{Automatic}           & \multicolumn{1}{c|}{0.394}           & \multicolumn{1}{c|}{334}        & \multicolumn{1}{c|}{114.8}   & \multicolumn{1}{c|}{$8.59 \time 10^{-3}$} \\ \hline
\end{tabular}
\caption{Results for compressed column for minimization of volume.}
\label{VCBtab}
\end{table}

\subsection{Cantilever results}\label{sec32}
The cantilever problem, Figure \ref{examples}b, is solved to minimize nonlinear end compliance, with a 40\% volume constraint, and a constraint on stability, where the minimum critical load factor is 2. The stability constraint is formulated using an eigenvalue analysis at a load factor of 1, which provides a good approximation of the true critical load factor, whilst being robust and efficient \citep{dunning2023stability}. The eigenvalue analysis also allows for multiple modes to be included via a K-S function to help avoid issues with mode switching. For all cantilever problems, 6 modes are used with as K-S aggregation parameter of 50. The problem is solved using the Method of Moving Asymptotes \citep{svanberg2002class}, with a conservative move limit of 0.05.

\begin{figure}[h]
    \centering
    \includegraphics[width=0.75\textwidth]{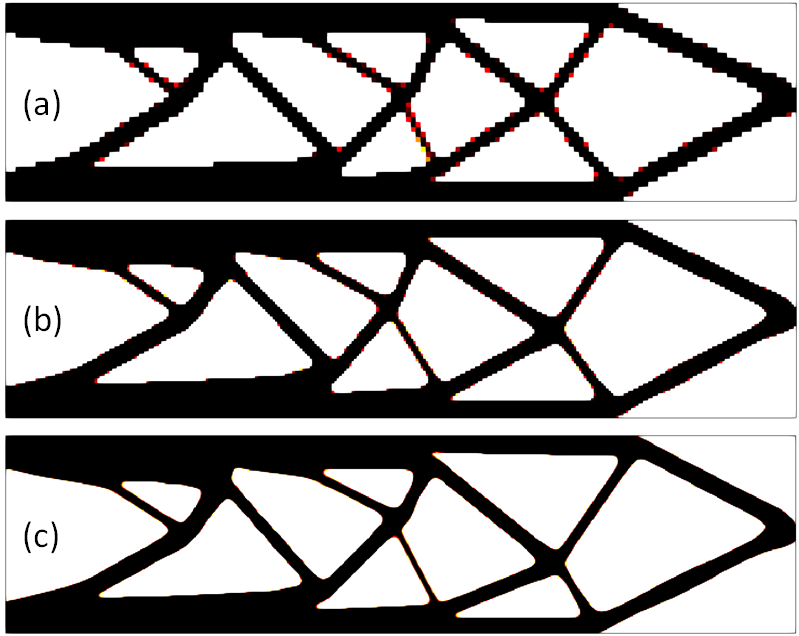}
    \caption{Cantilever results with load = $2 \times 10^{5}$ using mesh size: a) $160 \times 40$, b) $320 \times 80$, c) $640 \times 160$.}
    \label{Cant2Ksol}
\end{figure}

\begin{figure}[h]
    \centering
    \includegraphics[width=0.95\textwidth]{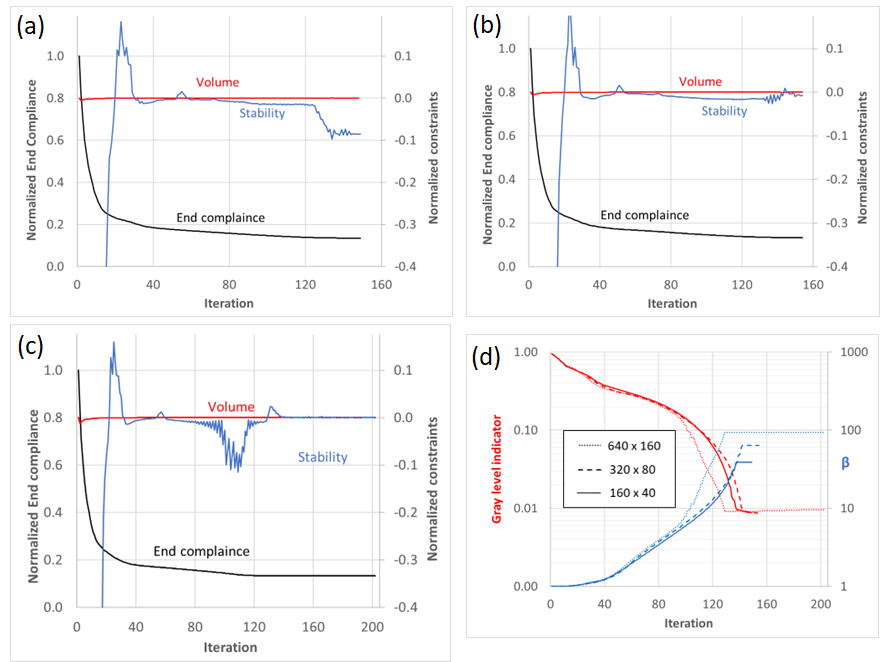}
    \caption{Cantilever with load = $2 \times 10^{5}$. Convergence history for mesh size: a) $160 \times 40$, b) $320 \times 80$, c) $640 \times 160$. d) Convergence of gray level indicator function, $G(x)$, and threshold projection parameter, $\beta$.}
    \label{2Kplot}
\end{figure}

The solutions for a load of $2 \times 10^{5}$ and the three different mesh sizes are shown in Figure \ref{Cant2Ksol}. The designs are very similar, due to the mesh independent filter, but with denser meshes producing smoother solutions. The automatic continuation scheme for $\beta$ is able to meet the gray level indicator stopping criterion, by choosing a larger $\beta$ value for a denser mesh - see Figure \ref{2Kplot}d and Table \ref{2Ktab}. It is also noticed that the value of $\beta$ follows a similar increase for about the first 90 iterations for all meshes. After this, $\beta$ is increased quicker and to a larger value for denser meshes to meet the gray level stopping criterion. Convergence histories for the three meshes are shown in Figure \ref{2Kplot}, which show smooth convergence, except for some oscillations in the stability constraint. This is mainly due to some mode switching, that could be improved by including more eigenvalues in the constraint formulation. However, feasible designs are obtained in a reasonable number of iterations. Finally, it is observed that the objective (end compliance) reduces as the mesh is refined, primarily due to the better resolution of topological features, but this comes at more computational cost, for both the finite elements analysis and number of iterations, see Table \ref{2Ktab}.

\begin{table}[h]
\begin{tabular}{|l|c|c|c|c|}
\hline
Mesh size        & End compliance & Iterations & $\beta$ & Gray value            \\ \hline
$160 \times 40$  & 57,068        & 149        & 38.7    & $8.77 \times 10^{-3}$ \\ \hline
$320 \times 80$  & 56,155        & 154        & 63.7    & $8.67 \times 10^{-3}$ \\ \hline
$640 \times 160$ & 55,840        & 202        & 93.2    & $9.59 \times 10^{-3}$ \\ \hline
\end{tabular}
\caption{Results for cantilever with load = $2 \times 10^{5}$.}
\label{2Ktab}
\end{table}

\begin{figure}[h]
    \centering
    \includegraphics[width=0.75\textwidth]{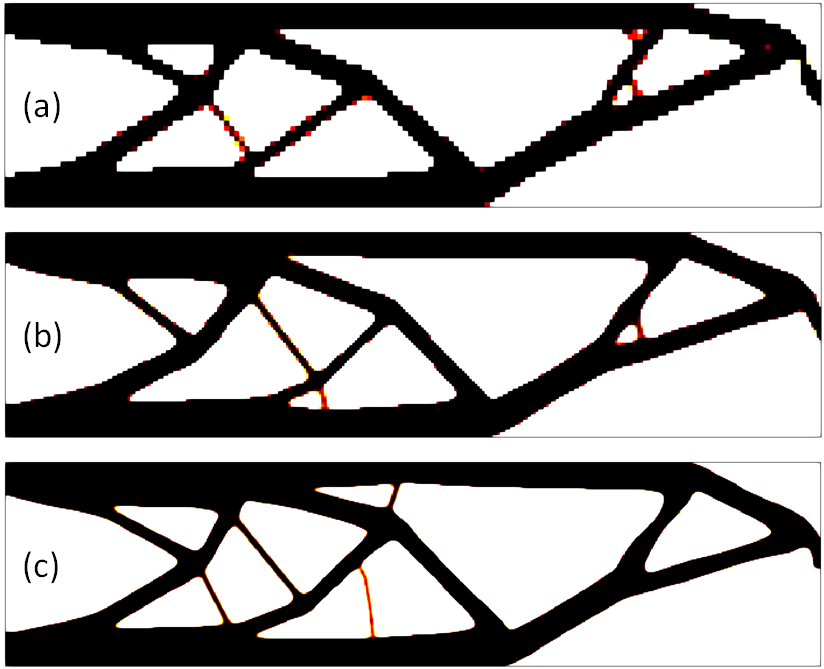}
    \caption{Cantilever results with load = $3 \times 10^{5}$ using mesh size: a) $160 \times 40$, b) $320 \times 80$, c) $640 \times 160$.}
    \label{Cant3Ksol}
\end{figure}

\begin{figure}[h]
    \centering
    \includegraphics[width=0.95\textwidth]{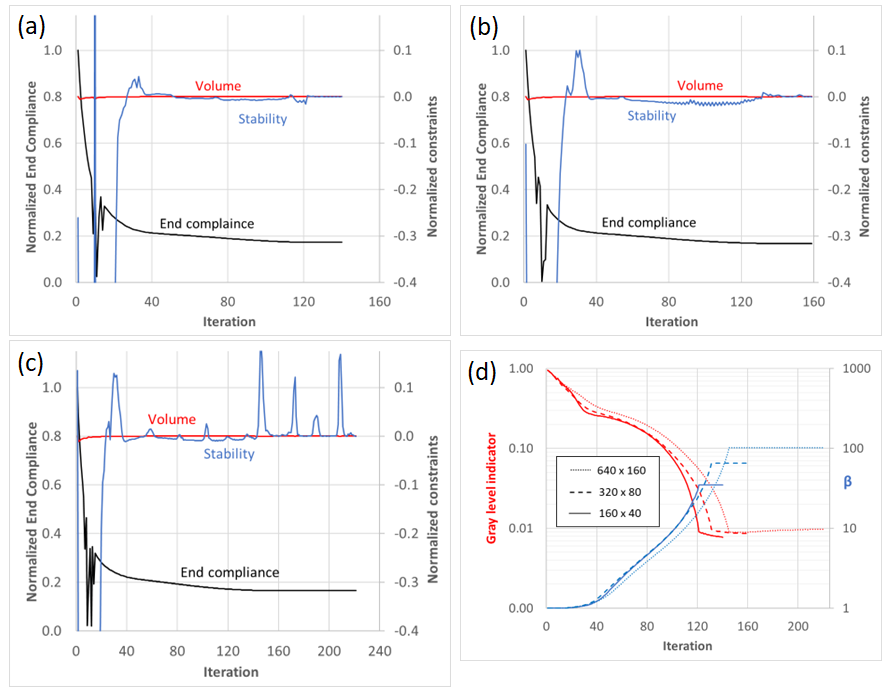}
    \caption{Cantilever with load = $3 \times 10^{5}$. Convergence history for mesh size: a) $160 \times 40$, b) $320 \times 80$, c) $640 \times 160$. d) Convergence of gray level indicator function, $G(x)$, and threshold projection parameter, $\beta$.}
    \label{3Kplot}
\end{figure}

Cantilever solutions for a load of $3 \times 10^{5}$ and the three different mesh sizes are shown in Figure \ref{Cant3Ksol}. Again, the designs are very similar, although the finest mesh ($640 \times 160$) has some additional thin members, which cannot be represented accurately using the coarser meshes. Compared with the results for a load of $2 \times 10^{5}$, similar observation are made about the final value of $\beta$ being larger for a finer mesh to meet the gray level stopping criterion, and that more iterations are needed for the finer mesh, but the objective function is lower - see Table \ref{3Ktab}. Convergence histories are also shown in Figure \ref{3Kplot}, which again show some oscillations in the stability constraints. This is most noticeable for the finest mesh ($640 \times 160$) and is caused by mode switching due to the longest thin member becoming even thinner. This can be eliminated by increasing the number of modes in the constraint formulation to 12, but the result for 6 modes is presented for consistency, and a feasible design is obtained anyway.

\begin{table}[h]
\begin{tabular}{|l|c|c|c|c|}
\hline
Mesh size        & End compliance & Iterations & $\beta$ & Gray value            \\ \hline
$160 \times 40$  & 134,415       & 140        & 35.0    & $7.73 \times 10^{-3}$ \\ \hline
$320 \times 80$  & 131,305       & 159        & 65.1    & $8.56 \times 10^{-3}$ \\ \hline
$640 \times 160$ & 129,422       & 221        & 101.5   & $9.67 \times 10^{-3}$ \\ \hline
\end{tabular}
\caption{Results for cantilever with load = $3 \times 10^{5}$.}
\label{3Ktab}
\end{table}

\section{Conclusions}\label{sec4}

A simple method is proposed to to achieve near binary designs for three-field density-based topology optimization, by automatically increasing the value of the threshold projection parameter during optimization. The method has two key parts. Firstly the increase in the threshold projection parameter is linked to the progress of the objective function, and secondly the gray level indicator function is used as a stopping criterion.

The method is tested on some 2D benchmark problems for linear buckling and a geometrically nonlinear case. The results show that the proposed method can achieve near binary solutions (as measured by the gray level indicator function). When compared with a conventional continuation approach, where the increase in $\beta$ is predetermined and not linked to optimization progress, improved results in terms of reduced grayness, and/or reduced number of iterations is observed. Further testing is required to determine the utility of the proposed method for a wider range of problems, including different physics, and when using different optimizers.

%%%%%%%%%%%%%%%

\bibliography{topopt}% common bib file

\end{document}